# A remark on a variable-coefficient Bernoulli equation based on auxiliary -equation method for nonlinear physical systems

*Zehra Pınar*[a] *Turgut Öziş*[b]

[a]*Namık Kemal University, Faculty of Arts and Science, Department of Mathematics, 59030Merkez-Tekirdağ, Turkey*

[b]*Ege University, Science Faculty, Department of Mathematics, 35100 Bornova-İzmir, Turkey*

Abstract: It is well recognized that in auxiliary equation methods, the exact solutions of different types of auxiliary equations may produce new types of exact travelling wave solutions to nonlinear partial differential equations in hand. In this study, we extend the class of auxiliary equations of classical Bernoulli equation which considered by various researchers [27, 31, 32, 33, 34, 35] to a variable-coefficient Bernoulli type equation. The proposed variable-coefficient Bernoulli type auxiliary equation produces many new solutions comparing to classical Bernoulli equation which produce two solutions only. Consequently, we introduce new exact travelling wave solutions of some physical systems in terms of these new solutions of the variable-coefficient Bernoulli type equation.

[a] Corresponding author: Zehra Pınar (zpinar@nku.edu.tr)

## 1. Introduction

The inspection of nonlinear wave phenomena of diffusion, convection, dispersion and dissipation appearing in physical systems is of great interest from both mathematical and physical points of view. In most case, the theoretical modelling based on nonlinear partial differential equations (NLPDEs) can accurately describe the wave dynamics of many physical systems. Of critical importance is to find closed form solutions for NLPDEs of physical significance. This could be a very complicated task and, in fact, is not always possible since in various realistic problems in physical systems. So, searching for some exact physically significant solutions is an important topic because of wide applications of NLPDEs in biology, chemistry, physics, fluid dynamics and other areas of engineering.

There are several theoretical results about local and global solutions of differential equations that establish existence, uniqueness etc., yet makes every effort to find exact solution formulas [1, 2].

Since many of the most useful techniques in analysis are formal or heuristic the trend in recent years has also been to justify and provide the new procedures or methods rigorously. Hence, over the past decades, a number of approximate methods for finding travelling wave solutions to nonlinear evolution equations have been proposed/or developed and furthermore modified. The solutions to various evolution equations have been found by one or other of these methods

Among all these methods, one of the prominent methods is so called auxiliary equation method.

The technique of this method consist of the solutions of the nonlinear evolution equations such that the target solutions of the nonlinear evolution equations can be expressed as a polynomial in an elementary function which satisfies a particular ordinary differential equation along with is named as auxiliary equation in general. Recently, to determine the solutions of nonlinear evolution equations, many exact solutions of various auxiliary equations have been utilized [23-28].

Recently, group of researchers [27, 31, 32, 33, 34, 35] have utilized the classical Bernoulli equation as an auxiliary equation to solve various nonlinear physical systems.

In this paper, we will examine the consequences of the choice of the auxiliary equation as variable coefficient Bernoulli type equation for determining the solutions of the nonlinear evolution equation in consideration and more we search for additional forms of new exact solutions of nonlinear differential equations which satisfying variable coefficient Bernoulli type equation(s).

## 2. Auxiliary equation technique and some remarks

The fundamental nature of the auxiliary equation technique is as follows. Let us have a nonlinear partial differential equation

$$P(u, u_x, u_t, u_{xx}, u_{xt}, u_{tt}, \cdots) = 0 \quad (1)$$

and let by means of an appropriate transformation this equation be reduced to nonlinear ordinary differential equation

$$Q\left(u, u_\xi, u_{\xi\xi}, u_{\xi\xi\xi}, \cdots\right) = 0. \quad (2)$$

For large class of the equations of the type (1) have exact solutions which can be constructed via finite series

$$u(\xi) = \sum_{i=0}^{N} a_i z^i(\xi) \quad (3)$$

Here, $a_i, (i = 0,1, \dots, N)$ are constants to be further determined, $N$ is an integer fixed by a balancing principle and elementary function $z(\xi)$ is the solution of some ordinary differential equation referred to as the auxiliary equation.

It is worth to point out that we happen to know the general solution(s), $z(\xi)$, of the auxiliary equation beforehand or we know at least exact analytical particular solutions of the auxiliary equation.

*The outline of the auxiliary equation method*:

A) Define the solution of Eq.(2) by the ansatz in form of finite series in Eq.(3) where $a_i, (i = 0,1, \dots, N)$ are constants to be further determined, $N$ is an integer fixed by a balancing principle and elementary function $z(\xi)$ is the solution of the auxiliary equation to be considered.

B) Substitute Eq.(3) into ordinary differential equation (2) to determine the coefficients $a_i, (i = 0,1, \dots, N)$ with the aid of symbolic computation.

C) Insert predetermined coefficients $a_i$ and elementary function $z(\xi)$ of the auxiliary equation into Eq.(3) to obtain travelling wave solutions of the nonlinear evolution equation in consideration.

It is very apparent that determining the elementary function $z(\xi)$ via auxiliary equation is crucial and plays very important role finding new travelling wave solutions of nonlinear evolution equations.

This fact, indeed, compel researchers to search for a novel auxiliary equations with exact solutions.

Recently, however, in most auxiliary equation methods, the elementary function $z(\xi)$ is defined as the solution(s) of the ordinary differential equations such as celebrity Bernoulli and Riccati equations which are well known nonlinear ordinary differential equations and their solutions can be expressed by elementary functions [27, 31, 32, 33, 34, 35].

In this study, we extend the classical Bernoulli equation [27, 31, 32, 33, 34, 35] to variable coefficient Bernoulli type equation(s)

$$\frac{dz}{d\xi} = P(\xi)z(\xi) + Q(\xi)\big(z(\xi)\big)^n \qquad (4)$$

where $P(\xi)$ and $Q(\xi)$ are polynomials, exponential functions and trigonometric functions respectively, $n$ is an integer and $n > 1$. Table 1 reports the new exact analytical solutions of the auxiliary equation (4) for distinct coefficient functions and for $n = 2$ under twenty distinct cases.

REMARK 1: It is remarkable to point out that the classical Bernoulli equation [27, 33] given below

$$\frac{dz}{d\xi} = az(\xi) + b[z(\xi)]^k \qquad (5)$$

where $k$ is an integer and $k > 1$ , give only a pair of solutions with respect to the constant coefficients $a$ and $b$ respectively. i.e.;

I) $z(\xi) = \sqrt[k-1]{\frac{aexp[a(k-1)(\xi+\xi_0)}{1-bexp[a(k-1)(\xi+\xi_0)}}$

for case a$< 0$ , $b > 0$ and

II) $z(\xi) = \sqrt[k-1]{-\frac{aexp[a(k-1)(\xi+\xi_0)}{1-bexp[a(k-1)(\xi+\xi_0)}}$

for $a > 0, b < 0$ where $\xi_0$ is a constant of integration[27, 33].

REMARK 2: In Table 1, the case 4 coincides with ref. [27] for $= A$ , $b = B$ and $C1$ is an integration constant and for $\xi_0 = 0$.

REMARK 3: In Table 1, the cases, expect the case 4, are novel solutions in the literature.

REMARK 4: In Table 1, some of the cases lead to solutions in complex plane. These complex solutions may have potential use for fluid dynamics, electromagnetism, electrical and electronics engineering, signal analysis, quantum mechanics etc. where for a convenient description for periodically varying solutions are needed.

## 3. Travelling Wave Solutions of b-equation in terms of variable coefficient Bernoulli Equation regarding as an ansatz.

Holm&Staley [36] studied the exchange of stability in the dynamics of solitary wave solutions under changes in the nonlinear balance in a $1+1$ evolutionary partial differential equation related both to shallow water waves and to turbulence. The family of b-equations are defined

$$\frac{\partial u}{\partial t} - \frac{\partial^3 u}{\partial x^2 \partial t} + (b+1)u\frac{\partial u}{\partial x} - b\frac{\partial u}{\partial x}\frac{\partial^2 u}{\partial x^2} - u\frac{\partial^3 u}{\partial x^3} = 0 \qquad (6)$$

where $u(x,t)$ denotes the velocity field.

For $b = 2$, Eq.(6) reduces to the Camassa–Holm equation and for $b = 3$, Eq.(6) reduces to the Degasperis–Processi equation [27]. The b-equations are of relevance for hydrodynamics [28–30] and are integrable only for the cases $b = 2$ (Camassa–Holm case) and $b = 3$ (Degasperis–Processi case).

Now, we consider the b-equation for $b = -2$,

$$\frac{\partial u}{\partial t} - \frac{\partial^3 u}{\partial x^2 \partial t} - u\frac{\partial u}{\partial x} + 2\frac{\partial u}{\partial x}\frac{\partial^2 u}{\partial x^2} - u\frac{\partial^3 u}{\partial x^3} = 0 \qquad (7)$$

To find the travelling wave solutions for Eq.(7), we use the wave variable $\xi = \mu(x - ct)$, where $c \neq 0$ and $\mu \neq 0$. The wave variable $\xi$ carries Eq. (7) into the ordinary differential equation

$$cu' - \mu^3 u''' - \mu u u' + 2\mu^3 u' u'' - \mu^3 u u''' = 0 \qquad (8)$$

From balancing principle and using (8), we have $N = 2$ , therefore the ansatz yields

$$u(\xi) = g_0 + g_1 z(\xi) + g_2 z^2(\xi) \qquad (9)$$

where $z(\xi)$ is a solution(s) of Eq. (4) ( see Table 1).

**Case 1:** In the event of $P(\xi) = (A\xi + B)^2$, $Q(\xi) = A\xi + B$ and $n = 2$, Bernoulli type equation is the following form

$$\frac{dz}{d\xi} = (A\xi + B)^2 z(\xi) + (A\xi + B\,)\big(z(\xi)\big)^2 \qquad (10)$$

where $A$ and $B$ are real constants. The solution of Eq.(10) is

$$z(\xi) = \frac{e^{(1/3 A^2 \xi^3 + A \xi^2 B + \xi B^2)}}{\int -e^{(1/3 A^2 \xi^3 + A \xi^2 B + \xi B^2)} (A\xi + B)\, d\xi + _C1}$$

and in the case of $A = 0$, the solution is reduced into the solution in [36]

$$z(\xi) = \frac{B}{-1 + _C1 B e^{-B^2 \xi}}$$

Hence, substituting Eqs.(9) and(10) into Eq.(8) and letting each coefficient of $z^i(\xi)\sqrt{a_2 z^2(\xi)+a_6 z^6(\xi)}$, $(0\le i\le 8)$ to be zero, we obtain

$$eq_6 := -14\,g_2 A - 22\,g_2 A^3\xi^3 - 66\,g_2 A^2\xi^2 B - 66\,g_2 A\xi B^2 - 22\,g_2 B^3 = 0$$

$$\begin{aligned} eq_5 := &-72\,g_0 g_2 A^2\xi^2 B - 72\,g_0 g_2 A\xi B^2 - 72\,c\,g_2 A^2\xi^2 B - 72\,c\,g_2 A\xi B^2 - 2\,g_2^{\,2} A\xi \\ &- 24\,g_0 g_2 B^3 - 24\,c\,g_2 B^3 - 2\,g_2^{\,2} B - 10\,g_2^{\,2} B^5 - 24\,g_0 g_2 A^3\xi^3 - 24\,c\,g_2 A^3\xi^3 \\ &- 50\,g_2^{\,2} A^4\xi^4 B - 100\,g_2^{\,2} A^3\xi^3 B^2 - 100\,g_2^{\,2} A^2\xi^2 B^3 - 50\,g_2^{\,2} A\xi B^4 - 12\,g_2^{\,2} A^3\xi^2 \\ &- 20\,g_2^{\,2}(A\xi+B)A^2\xi - 20\,g_2^{\,2}(A\xi+B)AB - 10\,g_2^{\,2} A^5\xi^5 - 24\,g_2^{\,2} A^2\xi B \\ &- 12\,g_2^{\,2} A B^2 = 0 \end{aligned}$$

$$\begin{aligned} eq_4 := &-216\,c\,g_2 A^3\xi^3 B - 324\,c\,g_2 A^2\xi^2 B^2 - 216\,c\,g_2 A\xi B^3 - 216\,g_0 g_2 A^3\xi^3 B \\ &- 324\,g_0 g_2 A^2\xi^2 B^2 - 216\,g_0 g_2 A\xi B^3 - 2\,g_2^{\,2} B^2 - 2\,g_2^{\,2} A^2\xi^2 - 54\,c\,g_2 B^4 \\ &- 54\,g_0 g_2 B^4 - 4\,g_2^{\,2} A\xi B - 54\,c\,g_2 A^4\xi^4 - 54\,g_0 g_2 A^4\xi^4 - 18\,g_0 g_2 AB \\ &- 18\,c\,g_2 AB - 16\,g_2^{\,2}(A\xi+B)AB^2 - 16\,g_2^{\,2}(A\xi+B)A^3\xi^2 - 18\,g_0 g_2 A^2\xi \\ &- 18\,c\,g_2 A^2\xi - 32\,g_2^{\,2}(A\xi+B)A^2\xi B - 4\,g_2^{\,2} A^2 = 0 \end{aligned}$$

$$\begin{aligned} eq_3 := &-38\,g_0 g_2 A^5\xi^5 - 38\,c\,g_2 A^5\xi^5 + 2\,c\,g_2 A\xi - 2\,g_0 g_2 A\xi - 38\,g_0 g_2 B^5 - 38\,c\,g_2 B^5 \\ &+ 2\,c\,g_2 B - 2\,g_0 g_2 B - 190\,g_0 g_2 A^4\xi^4 B - 380\,g_0 g_2 A^3\xi^3 B^2 - 380\,g_0 g_2 A^2\xi^2 B^3 \\ &- 190\,g_0 g_2 A\xi B^4 - 190\,c\,g_2 A^4\xi^4 B - 380\,c\,g_2 A^3\xi^3 B^2 - 380\,c\,g_2 A^2\xi^2 B^3 \\ &- 190\,c\,g_2 A\xi B^4 - 16\,g_0 g_2 AB^2 - 16\,c\,g_2 AB^2 - 32\,g_0 g_2 A^2\xi B - 32\,c\,g_2 A^2\xi B \\ &- 16\,g_0 g_2 A^3\xi^2 - 16\,c\,g_2 A^3\xi^2 - 28\,g_0 g_2(A\xi+B)AB - 28\,c\,g_2(A\xi+B)AB \\ &- 28\,g_0 g_2(A\xi+B)A^2\xi - 28\,c\,g_2(A\xi+B)A^2\xi = 0 \end{aligned}$$

$$\begin{aligned} eq_2 := &\,4\,c\,g_2 A\xi B - 48\,c\,g_2 A^5\xi^5 B - 120\,c\,g_2 A^4\xi^4 B^2 - 160\,c\,g_2 A^3\xi^3 B^3 \\ &- 120\,c\,g_2 A^2\xi^2 B^4 - 48\,c\,g_2 A\xi B^5 - 4\,g_0 g_2 A\xi B - 48\,g_0 g_2 A^5\xi^5 B \\ &- 120\,g_0 g_2 A^4\xi^4 B^2 - 160\,g_0 g_2 A^3\xi^3 B^3 - 120\,g_0 g_2 A^2\xi^2 B^4 - 48\,g_0 g_2 A\xi B^5 \\ &+ 2\,c\,g_2 B^2 - 8\,c\,g_2 B^6 - 2\,g_0 g_2 B^2 + 2\,c\,g_2 A^2\xi^2 - 2\,g_0 g_2 A^2\xi^2 - 8\,c\,g_2 A^6\xi^6 \\ &- 8\,g_0 g_2 B^6 - 8\,g_0 g_2 A^6\xi^6 - 48\,g_0 g_2(A\xi+B)A^2\xi B - 48\,c\,g_2(A\xi+B)A^2\xi B \\ &- 4\,g_0 g_2 A^2 - 4\,c\,g_2 A^2 - 24\,g_0 g_2(A\xi+B)A^3\xi^2 - 24\,c\,g_2(A\xi+B)A^3\xi^2 \\ &- 24\,g_0 g_2(A\xi+B)AB^2 - 24\,c\,g_2(A\xi+B)AB^2 = 0 \end{aligned}$$

$$eq_1 := c A^2 \xi^2 - c A^6 \xi^6 - g_0 A^2 \xi^2 - g_0 A^6 \xi^6 + 2 c A \xi B - 15 c A^4 \xi^4 B^2 - 20 c A^3 \xi^3 B^3 - 15 c A^2 \xi^2 B^4 - 6 c A \xi B^5 - 2 g_0 A \xi B - 15 g_0 A^4 \xi^4 B^2 - 20 g_0 A^3 \xi^3 B^3 - 15 g_0 A^2 \xi^2 B^4 - 6 g_0 A \xi B^5 - 6 c (A \xi + B) A B^2 - 6 g_0 (A \xi + B) A B^2 - 6 c (A \xi + B) A^3 \xi^2 - 6 g_0 (A \xi + B) A^3 \xi^2 - 12 c (A \xi + B) A^2 \xi B - 12 g_0 (A \xi + B) A^2 \xi B + c B^2 - c B^6 - g_0 B^2 - g_0 B^6 - 2 c A^2 - 2 g_0 A^2 - 6 c A^5 \xi^5 B - 6 g_0 A^5 \xi^5 B = 0$$

Solving the system by the aid of Maple 13, we can determine the coefficients:

$$g_2 = -(27 c B^4 - c A^2 \xi^2 + 4 c A^6 \xi^6 + g_0 A^2 \xi^2 + 4 g_0 A^6 \xi^6 - 2 c A \xi B + 60 c A^4 \xi^4 B^2 + 80 c A^3 \xi^3 B^3 + 60 c A^2 \xi^2 B^4 + 24 c A \xi B^5 + 2 g_0 A \xi B + 60 g_0 A^4 \xi^4 B^2 + 80 g_0 A^3 \xi^3 B^3 + 60 g_0 A^2 \xi^2 B^4 + 24 g_0 A \xi B^5 + 22 c A B^2 + 36 c A^2 B^2 \xi + 22 c A^3 \xi^2 + 36 c A^3 \xi^2 B + 44 c A^2 \xi B - c B^2 + 4 c B^6 + g_0 B^2 + 4 g_0 B^6 + 2 c A^2 + 2 g_0 A^2 + 24 c A^5 \xi^5 B + 24 g_0 A^5 \xi^5 B + 27 g_0 A^4 \xi^4 + 9 g_0 A B + 9 g_0 A^2 \xi + g_0 A \xi + 19 g_0 A^5 \xi^5 + 22 g_0 A B^2 + 22 g_0 A^3 \xi^2 + 12 c A B^3 + 12 c A^4 \xi^3 - c B + 36 g_0 A^2 B^2 \xi + 36 g_0 A^3 \xi^2 B + 27 c A^4 \xi^4 + 9 c A B + 9 c A^2 \xi + 19 c A^5 \xi^5 - c A \xi + 12 g_0 A B^3 + 12 g_0 A^4 \xi^3 + 108 c A^3 \xi^3 B + 162 c A^2 \xi^2 B^2 + 108 c A \xi B^3 + 95 c A^4 \xi^4 B + 190 c A^3 \xi^3 B^2 + 190 c A^2 \xi^2 B^3 + 95 c A \xi B^4 + 19 c B^5 + g_0 B + 19 g_0 B^5 + 27 g_0 B^4 + 108 g_0 A^3 \xi^3 B + 162 g_0 A^2 \xi^2 B^2 + 108 g_0 A \xi B^3 + 95 g_0 A^4 \xi^4 B + 190 g_0 A^3 \xi^3 B^2 + 190 g_0 A^2 \xi^2 B^3 + 95 g_0 A \xi B^4 + 44 g_0 A^2 \xi B) / (8 A B^3 + 2 A \xi B + 24 A^2 B^2 \xi + 24 A^3 \xi^2 B + B^2 + 2 A^2 + 8 A^4 \xi^3 + A^2 \xi^2)$$

$$A = A, B = B, g_1 = 0, g_0 = g_0$$

Substituting the above coefficients into ansatz (9) with the solution of Bernoulli type equation, we obtain one of new solution of b-equation. For the special values of parameters, the solution is given in Figure 1.

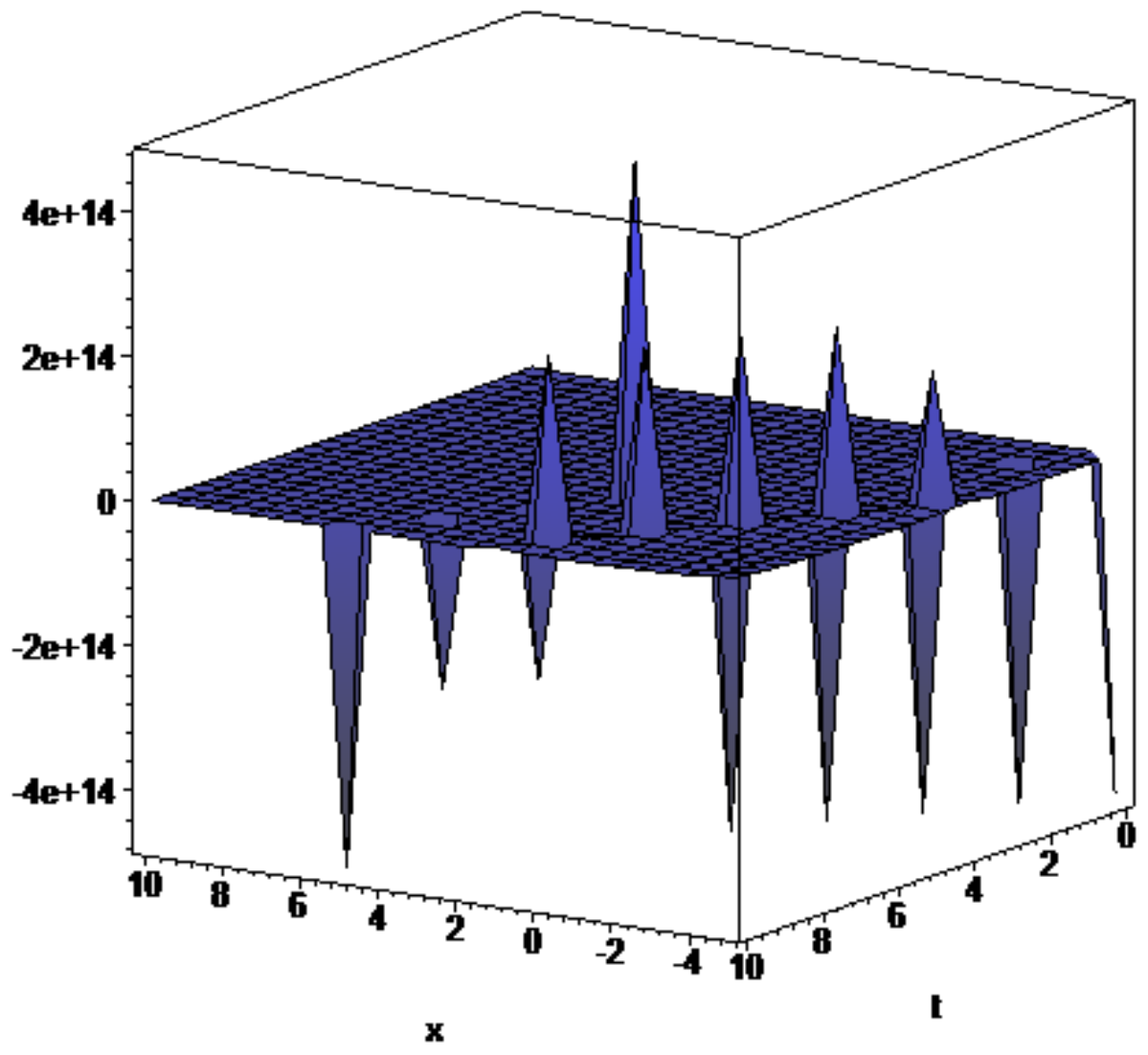

**Figure 1.** Graph of  solution of Eq.(8) for  $A = \frac{1}{4}, B = 1, c = \mu = 1, _C1 = 1.$

Now, we consider  $A = 0$, and our solution is reduced into  solution which is obtained by [27]  as seen in the Figure2.

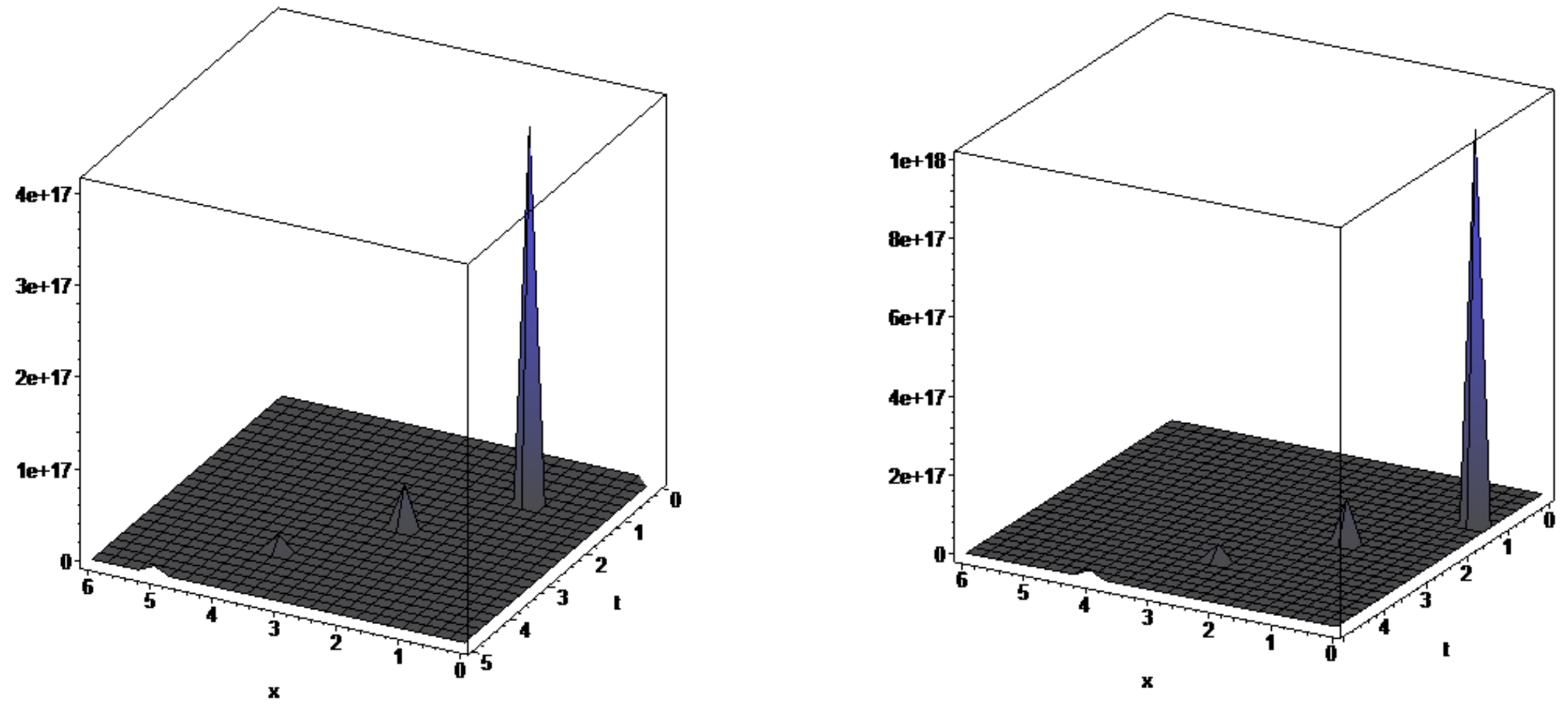

(a) (b)

**Figure 2.** (a) is the graph of reduced solution of Eq.(8) for  $A = 0, B = 1, c = 1,$ $\mu = 1, _C1 = 1$, (b) is the solution obtained by [27]

**Case 2:** In the event of $P(\xi) = (A\cos\xi)$, $Q(\xi) = (B\sin\xi)$ and $n = 2$, Bernoulli type equation is the following form

$$\frac{dz}{d\xi} = (A\cos\xi)z(\xi) + (B\sin\xi)\big(z(\xi)\big)^2 \qquad (11)$$

where $A$ and $B$ are real constants. The solution of Eq.(11) is

$$z(\xi) = \frac{e^{(A\sin(\xi))}}{\int -e^{(A\sin(\xi))} B\sin(\xi)\, d\xi + _C1}$$

Hence, substituting Eqs.(11) and (13) into Eq.(10) and letting each coefficient of $z^i(\xi)\sqrt{a_2 z^2(\xi) + a_6 z^6(\xi)}$, $(0 \le i \le 8)$ to be zero, we obtain

$$eq_{61} := -22\, g_2 A\cos(\xi)\sin(\xi) - 20\, g_1 B\sin(\xi)^2 - 14\, g_2\cos(\xi) = 0$$

$$\begin{aligned} eq_{51} := {} & -24\, g_0 g_2 B^2\sin(\xi)^3 - 40\, g_1 A\cos(\xi) g_2 B\sin(\xi)^2 - 24\, c\, g_2 B^2\sin(\xi)^3 \\ & - 20\, g_1 B\sin(\xi) g_2\cos(\xi) - 12\, g_2^2 A\cos(\xi)^2 + 10\, g_2^2\sin(\xi)^2 A \\ & - 10\, g_2^2 A^2\cos(\xi)^2\sin(\xi) - 4\, g_1^2 B^2\sin(\xi)^3 = 0 \end{aligned}$$

$$\begin{aligned} eq_4 := {} & -54\, g_0 g_2 A\cos(\xi) B^2\sin(\xi)^2 - 23\, g_1 A^2\cos(\xi)^2 g_2 B\sin(\xi) \\ & - 54\, c\, g_2 A\cos(\xi) B^2\sin(\xi)^2 - 18\, g_0 g_2 B^2\sin(\xi)\cos(\xi) + 14\, g_1 B\sin(\xi)^2 g_2 A \\ & + 8\, g_2^2 A^2\cos(\xi)\sin(\xi) - 6\, g_0 g_1 B^3\sin(\xi)^3 - 17\, g_1 A\cos(\xi)^2 g_2 B \\ & - 18\, c\, g_2 B^2\sin(\xi)\cos(\xi) - 5\, g_1^2 B^2\sin(\xi)\cos(\xi) - 7\, g_1^2 A\cos(\xi) B^2\sin(\xi)^2 \\ & - 6\, c\, g_1 B^3\sin(\xi)^3 = 0 \end{aligned}$$

$$\begin{aligned} eq_3 := {} & 4\, c\, g_2 B\sin(\xi) - 16\, g_0 g_2 A\cos(\xi)^2 B - 38\, g_0 g_2 A^2\cos(\xi)^2 B\sin(\xi) \\ & - 12\, g_0 g_1 A\cos(\xi) B^2\sin(\xi)^2 + 14\, g_0 g_2 B\sin(\xi)^2 A - 38\, c\, g_2 A^2\cos(\xi)^2 B\sin(\xi) \\ & - 16\, c\, g_2 A\cos(\xi)^2 B - 6\, g_0 g_1 B^2\cos(\xi)\sin(\xi) - 3\, g_1^2 A^2\cos(\xi)^2 B\sin(\xi) \\ & - 6\, c\, g_1 B^2\cos(\xi)\sin(\xi) + 11\, g_1 A^2\cos(\xi) g_2\sin(\xi) - 3\, g_1 A^3\cos(\xi)^3 g_2 \\ & + 3\, g_1^2 B\sin(\xi)^2 A - 4\, g_1^2 A\cos(\xi)^2 B - 12\, c\, g_1 A\cos(\xi) B^2\sin(\xi)^2 \\ & + 14\, c\, g_2 B\sin(\xi)^2 A = 0 \end{aligned}$$

$$\begin{aligned} eq_2 := {} & 2\, c\, g_1 B\sin(\xi) + 4\, c\, g_1 A\sin(\xi)^2 B - 7\, g_0 g_1 A^2\cos(\xi)^2 B\sin(\xi) + 4\, c\, g_2 A\cos(\xi) \\ & + 12\, g_0 g_2 A^2\cos(\xi)\sin(\xi) - 5\, g_0 g_1 A\cos(\xi)^2 B - 8\, c\, g_2 A^3\cos(\xi)^3 \\ & + 4\, g_0 g_1 A\sin(\xi)^2 B + 2\, g_1^2 A^2\cos(\xi)\sin(\xi) - 7\, c\, g_1 A^2\cos(\xi)^2 B\sin(\xi) \\ & + 12\, c\, g_2 A^2\cos(\xi)\sin(\xi) - 5\, c\, g_1 A\cos(\xi)^2 B - 8\, g_0 g_2 A^3\cos(\xi)^3 = 0 \end{aligned}$$

$$eq_1 := c\, g_1 - c\, A^2\cos(\xi)^2 + 3\, g_0 A\sin(\xi) + 3\, c\, A\sin(\xi) + c - g_0 A^2\cos(\xi)^2 = 0$$

Solving the system by the aid of Maple 13, the coefficients are determined

$$g_1 := 0$$

$$A = -\frac{7}{11}\frac{1}{\sin(\xi)}$$

$$B = -\frac{7}{132}(-1720\cos(\xi)\sin(\xi)^4 + 2325\cos(\xi) - 3286\cos(\xi)^3 + ($$
$$2958400\cos(\xi)^2\sin(\xi)^8 - 119203584\cos(\xi)^2\sin(\xi)^4 + 11303840\cos(\xi)^4\sin(\xi)^4$$
$$+ 83698665\cos(\xi)^2 - 71920092\cos(\xi)^4 + 10797796\cos(\xi)^6$$
$$+ 19475904\sin(\xi)^{10} - 45972192\sin(\xi)^8 + 72468480\sin(\xi)^4 - 19475904\sin(\xi)^6$$
$$- 26496288 + 32912544\cos(\xi)^2\sin(\xi)^6)^{(1/2)}) \big/ ($$
$$(674\sin(\xi)^4 - 674 + 1139\cos(\xi)^2)\sin(\xi)^2)$$

$$g_2 := -\frac{121}{6} v\sin(\xi)^2(78832\cos(\xi)\sin(\xi)^6 - 66732\cos(\xi)\sin(\xi)^2$$
$$+ 125632\cos(\xi)^3\sin(\xi)^2 - 59912\cos(\xi)\sin(\xi)^4 + 41157\cos(\xi) - 89486\cos(\xi)^3$$
$$- 31(2958400\cos(\xi)^2\sin(\xi)^8 - 119203584\cos(\xi)^2\sin(\xi)^4$$
$$+ 11303840\cos(\xi)^4\sin(\xi)^4 + 83698665\cos(\xi)^2 - 71920092\cos(\xi)^4$$
$$+ 10797796\cos(\xi)^6 + 19475904\sin(\xi)^{10} - 45972192\sin(\xi)^8 + 72468480\sin(\xi)^4$$
$$- 19475904\sin(\xi)^6 - 26496288 + 32912544\cos(\xi)^2\sin(\xi)^6)^{(1/2)} + 20\sin(\xi)^2($$
$$2958400\cos(\xi)^2\sin(\xi)^8 - 119203584\cos(\xi)^2\sin(\xi)^4 + 11303840\cos(\xi)^4\sin(\xi)^4$$
$$+ 83698665\cos(\xi)^2 - 71920092\cos(\xi)^4 + 10797796\cos(\xi)^6$$
$$+ 19475904\sin(\xi)^{10} - 45972192\sin(\xi)^8 + 72468480\sin(\xi)^4 - 19475904\sin(\xi)^6$$
$$- 26496288 + 32912544\cos(\xi)^2\sin(\xi)^6)^{(1/2)}) \big/ (\cos(\xi)(-14005046\sin(\xi)^4$$
$$+ 2193870 - 3707445\cos(\xi)^2 - 4968728\sin(\xi)^6 + 4968728\sin(\xi)^2$$
$$- 8396708\cos(\xi)^2\sin(\xi)^2 + 11811176\sin(\xi)^8 + 19959836\cos(\xi)^2\sin(\xi)^4))$$

$$g_0 := -\frac{1}{12}(-22785 + 4661\sin(\xi)^2 + 41114\sin(\xi)^4) v(78832\cos(\xi)\sin(\xi)^6$$
$$- 66732\cos(\xi)\sin(\xi)^2 + 125632\cos(\xi)^3\sin(\xi)^2 - 59912\cos(\xi)\sin(\xi)^4$$
$$+ 41157\cos(\xi) - 89486\cos(\xi)^3 - 31(2958400\cos(\xi)^2\sin(\xi)^8$$
$$- 119203584\cos(\xi)^2\sin(\xi)^4 + 11303840\cos(\xi)^4\sin(\xi)^4 + 83698665\cos(\xi)^2$$
$$- 71920092\cos(\xi)^4 + 10797796\cos(\xi)^6 + 19475904\sin(\xi)^{10} - 45972192\sin(\xi)^8$$
$$+ 72468480\sin(\xi)^4 - 19475904\sin(\xi)^6 - 26496288 + 32912544\cos(\xi)^2\sin(\xi)^6)\hat{}$$
$$(1/2) + 20\sin(\xi)^2(2958400\cos(\xi)^2\sin(\xi)^8 - 119203584\cos(\xi)^2\sin(\xi)^4$$
$$+ 11303840\cos(\xi)^4\sin(\xi)^4 + 83698665\cos(\xi)^2 - 71920092\cos(\xi)^4$$
$$+ 10797796\cos(\xi)^6 + 19475904\sin(\xi)^{10} - 45972192\sin(\xi)^8 + 72468480\sin(\xi)^4$$
$$- 19475904\sin(\xi)^6 - 26496288 + 32912544\cos(\xi)^2\sin(\xi)^6)^{(1/2)}) \big/ \Big(($$
$$-14005046\sin(\xi)^4 + 2193870 - 3707445\cos(\xi)^2 - 4968728\sin(\xi)^6$$
$$+ 4968728\sin(\xi)^2 - 8396708\cos(\xi)^2\sin(\xi)^2 + 11811176\sin(\xi)^8$$
$$+ 19959836\cos(\xi)^2\sin(\xi)^4)\Big(98\cos(\xi)\sin(\xi)^2 - 98\cos(\xi) + \frac{217}{12}($$

$$1720\cos(\xi)\sin(\xi)^4 - 2325\cos(\xi) + 3286\cos(\xi)^3 - (2958400\cos(\xi)^2\sin(\xi)^8$$
$$- 119203584\cos(\xi)^2\sin(\xi)^4 + 11303840\cos(\xi)^4\sin(\xi)^4 + 83698665\cos(\xi)^2$$
$$- 71920092\cos(\xi)^4 + 10797796\cos(\xi)^6 + 19475904\sin(\xi)^{10} - 45972192\sin(\xi)^8$$
$$+ 72468480\sin(\xi)^4 - 19475904\sin(\xi)^6 - 26496288 + 32912544\cos(\xi)^2\sin(\xi)^6)\char`^$$
$$^{(1/2)}) \Big/ (674\sin(\xi)^4 - 674 + 1139\cos(\xi)^2) - \frac{35}{3}(1720\cos(\xi)\sin(\xi)^4$$
$$- 2325\cos(\xi) + 3286\cos(\xi)^3 - (2958400\cos(\xi)^2\sin(\xi)^8$$
$$- 119203584\cos(\xi)^2\sin(\xi)^4 + 11303840\cos(\xi)^4\sin(\xi)^4 + 83698665\cos(\xi)^2$$
$$- 71920092\cos(\xi)^4 + 10797796\cos(\xi)^6 + 19475904\sin(\xi)^{10} - 45972192\sin(\xi)^8$$
$$+ 72468480\sin(\xi)^4 - 19475904\sin(\xi)^6 - 26496288 + 32912544\cos(\xi)^2\sin(\xi)^6)\char`^$$
$$^{(1/2)})\sin(\xi)^2 \Big/ (674\sin(\xi)^4 - 674 + 1139\cos(\xi)^2)\Big)\Big)$$

Substituting the above coefficients into ansatz (11) with the solution of Bernoulli type equation, we obtain one of new solution of b-equation. For the special values of parameters, the solution is shown in Figure 3.

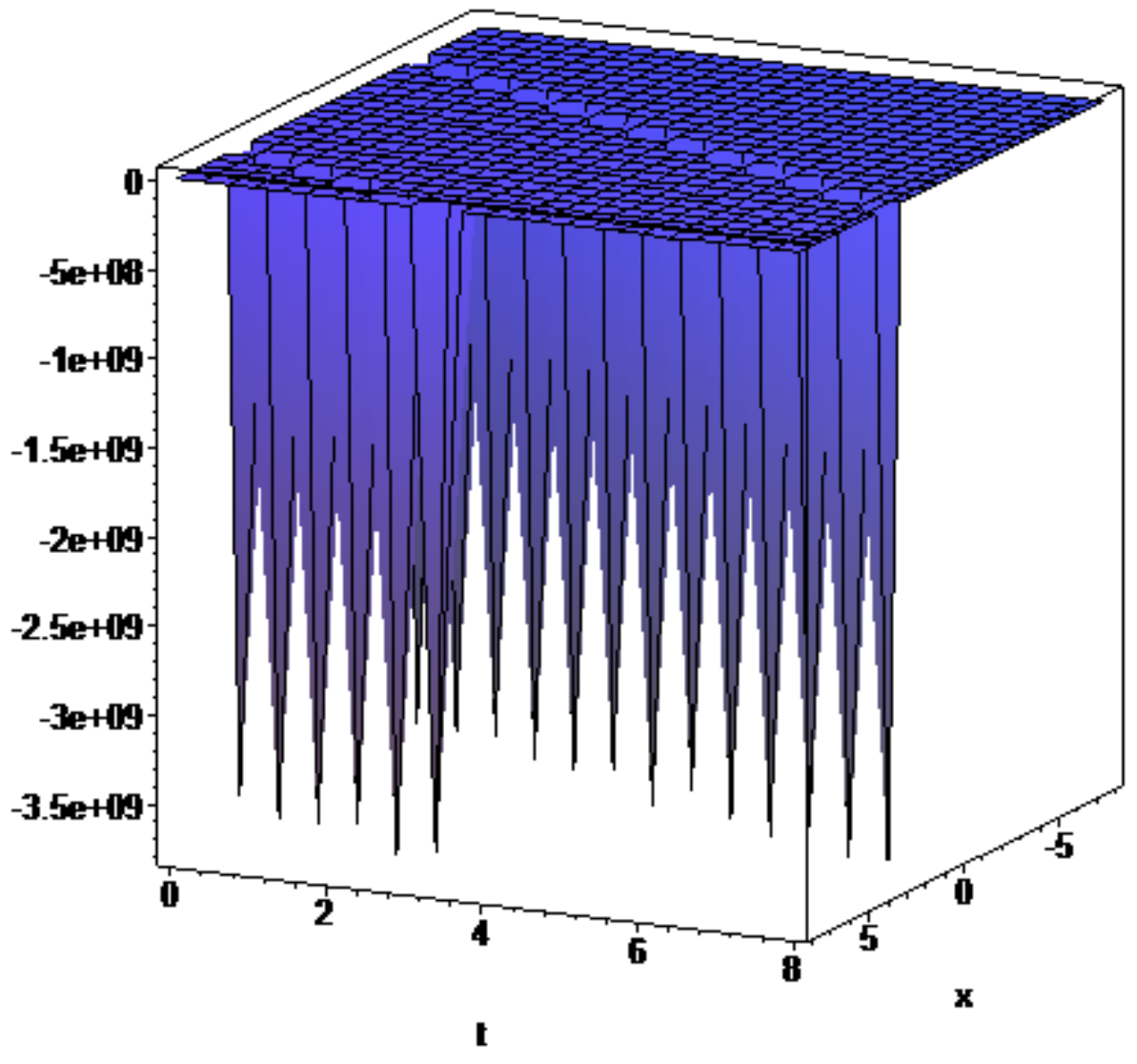

**Figure 3.** Graph of solution of Eq.(8) for $c = \mu = 1, _C1 = 1$.

Due to space limitation on the manuscript we only introduced two distinct cases to remark the aplicability of variable coefficient Berneoulli equation regarding as new ansatz to obtain new travelling solutions of nonlinear physical systems. In a similar manner, other solutions ( see Table 1) can also be used but we left it as an interesting implementation to the readers.

## 4. Conclusion

As is seen, the key idea of obtaining new travelling wave solutions for the nonlinear equations is using different type auxiliary equation. In this letter, the exact solutions of the variable coefficient Bernoulli type equation where $n$ is an integer and $n > 1$, $P(\xi)$ and $Q(\xi)$ are real functions (Eq. (4)) i.e.

$$\frac{dz}{d\xi} = P(\xi)z(\xi) + Q(\xi)\big(z(\xi)\big)^n$$

is utilized. Using the solutions of variable coefficient Bernoulli equation, we have successfully obtained a number of new exact periodic solutions of the b-equation by employing the solutions of the variable coefficient Bernoulli type equation regarding as an auxiliary equation in proposed method.

For comparative purposes, it is noticeable that Vitanov *et al* in [27, 31, 32, 33, 34] have utilized Eq.(4) for $P(\xi)$ and $Q(\xi)$ are constants, therefore their solutions are limited in the quantity to the order considered. In this letter, we have obtained solutions of the variable coefficient Bernoulli equation (Eq.(4)) as an auxiliary equation for distinct cases to obtain new travelling solutions of the nonlinear equation in hand. For that reason, our approach give additional new solutions beside the solutions obtained by using classical constant coefficient Bernoulli equation which already used aforementioned authors.

However, it is well known that different types of auxiliary equations produce new travelling wave solutions for many nonlinear problems. The presented method could lead to finding new exact travelling wave solutions for other nonlinear problems.

**Table 1.** The solutions of the variable coefficient Bernoulli type equation (4) regarding as an auxiliary equation

| Cases | $P(\xi)$ | $Q(\xi)$ | $z(\xi)$ |
|---|---|---|---|
| **Case 1** | $(A\xi + B)^2$ | $A\xi + B$ | $z(\xi) = \dfrac{e^{(1/3 A^2 \xi^3 + A \xi^2 B + \xi B^2)}}{\int -e^{(1/3 A^2 \xi^3 + A \xi^2 B + \xi B^2)} (A \xi + B)\, d\xi + _C1}$ |
| **Case 2** | $A \cos \xi$ | $B \sin \xi$ | $z(\xi) = \dfrac{e^{(A \sin(\xi))}}{\int -e^{(A \sin(\xi))} B \sin(\xi)\, d\xi + _C1}$ |
| **Case 3** | $A$ | $(C\xi + B)^2$ | $z(\xi) = \dfrac{A^3}{-A^2 B^2 + 2 A C B - 2 C^2 - 2 A^2 C \xi B + 2 A C^2 \xi - C^2 \xi^2 A^2 + e^{(-A \xi)}\, _C1\, A^3}$ |
| **Case 4** | $A$ | $B$ | $z(\xi) = \dfrac{A}{-B + e^{(-A \xi)}\, _C1\, A}$ |
| **Case 5** | $A$ | $e^{C\xi}$ | $z(\xi) = \dfrac{A + C}{-e^{(C \xi)} + e^{(-A \xi)}\, _C1\, A + e^{(-A \xi)}\, _C1\, C}$ |
| **Case 6** | $A$ | $e^{C\xi} + B$ | $z(\xi) = \dfrac{(A + C) A}{-A e^{(C \xi)} - A B - B C + e^{(-A \xi)}\, _C1\, A^2 + e^{(-A \xi)}\, _C1\, A C}$ |
| **Case 7** | $A \sin \xi$ | $C\xi + B$ | $z(\xi) = \dfrac{e^{(-A \cos(\xi))}}{\int -e^{(-A \cos(\xi))} (C \xi + B)\, d\xi + _C1}$ |
| **Case 8** | $A \cos \xi$ | $C\xi + B$ | $z(\xi) = \dfrac{e^{(A \sin(\xi))}}{\int -e^{(A \sin(\xi))} (C \xi + B)\, d\xi + _C1}$ |
| **Case 9** | $e^{C\xi}$ | $A$ | $z(\xi) = \dfrac{C e^{\left(\frac{e^{(C \xi)}}{C}\right)}}{\mathrm{Ei}\left(1, -\frac{e^{(C \xi)}}{C}\right) A + _C1\, C}$ |
| **Case10** | $e^{C\xi}$ | $A\xi + B$ | $z(\xi) = \dfrac{e^{\left(\frac{e^{(C \xi)}}{C}\right)}}{\int -e^{\left(\frac{e^{(C \xi)}}{C}\right)} (A \xi + B)\, d\xi + _C1}$ |

| | | | | | |
|---|---|---|---|---|---|
| **Case11** | | | $A\xi + B$ | $e^{C\xi}$ | $z(\xi) = \frac{\sqrt{-2A}\, e^{(1/2 A \xi^2 + B\xi)}}{\sqrt{\pi}\, e^{\left(-\frac{(B+C)^2}{2A}\right)} \operatorname{erf}\left(\frac{A\xi + B + C}{\sqrt{-2A}}\right) + _C1 \sqrt{-2A}}$ |
| **Case12** | | $A \sin \xi$ | | $e^{C\xi}$ | $z(\xi) = \frac{e^{(-A\cos(\xi))}}{\int -e^{(-A\cos(\xi))}\, e^{(C\xi)}\, d\xi + _C1}$ |
| **Case13** | $e^{C\xi}$ | | $A \sin \xi$ | | $z(\xi) = \frac{e^{\left(\frac{e^{(C\xi)}}{C}\right)}}{\int -e^{\left(\frac{e^{(C\xi)}}{C}\right)} \sin(\xi)\, A\, d\xi + _C1}$ |
| **Case14** | $A \cos \xi$ | | $e^{C\xi}$ | | $z(\xi) = \frac{e^{(A\sin(\xi))}}{\int -e^{(A\sin(\xi))}\, e^{(C\xi)}\, d\xi + _C1}$ |
| **Case15** | $e^{C\xi}$ | | $A \cos \xi$ | | $z(\xi) = \frac{e^{\left(\frac{e^{(C\xi)}}{C}\right)}}{\int -e^{\left(\frac{e^{(C\xi)}}{C}\right)} \cos(\xi)\, A\, d\xi + _C1}$ |
| **Case16** | $A \sin \xi$ | | $B \cos \xi$ | | $z(\xi) = \frac{e^{(-A\cos(\xi))}}{\int -e^{(-A\cos(\xi))}\, B\cos(\xi)\, d\xi + _C1}$ |
| **Case17** | | $(C\xi + B)^2$ | | $A$ | $z(\xi) = \frac{\sqrt{-2C}\, e^{(1/2 C \xi^2 + B\xi)}}{A\sqrt{\pi}\, e^{\left(-\frac{B^2}{2C}\right)} \operatorname{erf}\left(\frac{C\xi + B}{\sqrt{-2C}}\right) + _C1 \sqrt{-2C}}$ |
| **Case18** | | $e^{C\xi}$ | | $e^{B\xi}$ | $z(\xi) = \frac{e^{\left(\frac{e^{(C\xi)}}{C}\right)}}{\int -e^{\left(\frac{e^{(C\xi)}}{C}\right)} e^{(B\xi)}\, d\xi + _C1}$ |
| **Case19** | $C\xi + B$ | | $A \cos \xi$ | | $z(\xi) = \frac{2\sqrt{-2C}\, e^{\left(\frac{1}{2}C\xi^2 + B\xi\right)}}{A\sqrt{\pi}\, e^{\left(-\frac{(B+I)^2}{2C}\right)} \operatorname{erf}\left(\frac{C\xi + B + I}{\sqrt{-2C}}\right) + A\sqrt{\pi}\, e^{\left(-\frac{(B-I)^2}{2C}\right)} \operatorname{erf}\left(\frac{C\xi + B - I}{\sqrt{-2C}}\right) + 2_C1\sqrt{-2C}}$ |
| **Case20** | $C\xi + B$ | | $A \sin \xi$ | | $z(\xi) = -\frac{2\sqrt{-2C}\, e^{\left(\frac{1}{2}C\xi^2 + B\xi\right)}}{A\sqrt{\pi}\, e^{\left(-\frac{(B+I)^2}{2C}\right)} \operatorname{erf}\left(\frac{C\xi + B + I}{\sqrt{-2C}}\right) I - A\sqrt{\pi}\, e^{\left(-\frac{(B-I)^2}{2C}\right)} \operatorname{erf}\left(\frac{C\xi + B - I}{\sqrt{-2C}}\right) I + 2_C1\sqrt{-2C}}$ |

*The solutions of Bernoulli equation are obtained for $n = 2$.

**List of figures:**

**List of tables:**

**Table 1**. The solutions of the variable coefficient Bernoulli type equation (4) regarding as an auxiliary equation